\documentclass[11pt,a4paper,leqno]{amsart}

\usepackage[latin1]{inputenc}
\usepackage[T1]{fontenc}
\usepackage{amsfonts}
\usepackage{amsmath}
\usepackage{amssymb}
\usepackage{eurosym}
\usepackage{mathrsfs}
\usepackage{palatino}
\usepackage{color}
\usepackage{xcolor}
\usepackage{esint}
\usepackage{url}
\usepackage{hyperref}

\newcommand{\R}{\mathbb{R}}

\newcommand{\C}{\mathbb{C}}

\newcommand{\No}{\mathbb{N}_0}
\newcommand{\Z}{\mathbb{Z}}

\numberwithin{equation}{section}

\newcommand{\ud}[0]{\,\mathrm{d}}

\newcommand{\abs}[1]{|#1|}
\newcommand{\babs}[1]{\big|#1\big|}
\newcommand{\Babs}[1]{\Big|#1\Big|}
\newcommand{\norm}[2]{\|#1\|_{#2}}

\newcommand{\Norm}[2]{\|#1\|_{#2}}

\newcommand{\BNorm}[2]{\Big\|#1\Big\|_{#2}}

\newcommand{\ave}[1]{\langle #1\rangle}

\newcommand{\BMO}[0]{\operatorname{BMO}}
\newcommand{\supp}[0]{\operatorname{supp}}

\newcommand{\Exp}[0]{\mathbb{E}}

\theoremstyle{plain}
\newtheorem{thm}{Theorem}[section]
\newtheorem{lem}[thm]{Lemma}
\newtheorem{prop}[thm]{Proposition}
\newtheorem{cor}[thm]{Corollary}

\theoremstyle{definition}
\newtheorem{defn}[thm]{Definition}

\theoremstyle{remark}
\newtheorem{rem}[thm]{Remark}

\title{A sparse bilinear form with mean oscillations}

\author{Li Jinsong}

\address{Center for Applied Mathematics, Tianjin University, Weijin Road 92, 300072 Tianjin, China}
\email{ljs@tju.edu.cn}

\makeatletter
\@namedef{subjclassname@2020}{%
  \textup{2020} Mathematics Subject Classification}
\makeatother


\begin{document}

\begin{abstract}
In this article, we establish an oscillatory sparse domination for the bilinear forms associated with Calder\'on-Zygmund operators. We further characterize this oscillatory sparse domination in terms of either the boundedness of $T1$ or local endpoint Orlicz norm estimates. Our method also applies to recover Lerner's Sobolev-type inequality in \cite{L2026}.
\end{abstract}
\maketitle
\section{Introduction}
Sparse domination has become a central tool in modern harmonic analysis since Lerner \cite{L2013} provided an alternative proof of the $A_2$ theorem, which was based on norm domination of Calder\'on-Zygmund operators by positive sparse operators. Later, this norm estimate was strengthened to pointwise sparse domination in the works of Conde-Alonso and Rey \cite{CJR2016}, Lerner \cite{L2016}, and Lerner and Nazarov \cite{LN2019}. More precisely, if $T$ is a Calder\'on-Zygmund operator, then for every compactly supported function $f\in L^1(\R^d)$ there exists an $\eta$-sparse family $\mathcal S$ such that
\begin{equation}\label{1.1}
\abs{Tf(x)}
\lesssim
\sum_{Q\in\mathcal S}\ave{\abs{f}}_Q 1_Q(x), \,\text{ where } \ave{f}_Q:=\abs{Q}^{-1}\int_Qf,
\end{equation}
and the $\eta$-sparseness means that for each $Q\in\mathcal{S}$, there exist pairwise disjoint measurable subsets $E_Q\subset Q$ with $\abs{E_Q}\geq\eta\abs{Q}$. In all sparse bounds below, the sparsity parameter and the implicit constants are independent of the test functions.

Subsequently, sparse domination was found to be applicable to various other classes of operators. However, the natural sparse control is often no longer a pointwise estimate for $Tf$, but rather a bilinear-form estimate for the pairing $\langle Tf,g\rangle$. For instance, Conde-Alonso, Culiuc, Di Plinio, and Ou \cite{CCDO2017} established a sparse bilinear form for rough homogeneous singular integrals and Bochner-Riesz means at the critical index, with the usual $L^1$-averages replaced by suitable $L^r$-averages. Sparse bilinear form for oscillatory singular integrals and their maximal truncations were obtained by Lacey and Spencer \cite{LS2017} and by Krause and Lacey \cite{KL2020}, respectively. These results demonstrate the flexibility of bilinear sparse forms in settings where classical pointwise domination is unavailable. 

Recently, Conde-Alonso, Lorist, and Rey \cite{CLR2026} developed a cancellative sparse domination principle in which the usual size averages are replaced by local quantities that retain information about the cancellation of the input function. A related, but distinct, refinement was subsequently obtained by Lerner \cite{L2026}. Building on an earlier result of Benea and Bernicot \cite{BB2018}, he proved that if $T$ is a Calder\'on-Zygmund operator satisfying $T(1)=0$, then 
\begin{equation}\label{1.2}
\abs{Tf(x)}
\lesssim
\sum_{Q\in\mathcal S}\Omega(f;Q) 1_Q(x)\quad\text{ where }\,\Omega(f;Q):=\ave{\abs{f-\ave{f}_Q}}_Q.
\end{equation}
Such pointwise sparse domination immediately yields a corresponding sparse bilinear form, namely
\[
\abs{\ave{Tf_1,f_2}}\lesssim
\sum_{Q\in\mathcal S}\Omega(f_1;Q) \ave{\abs{f_2}}_Q\abs{Q}.
\]

In contemporaneous work, Di Plinio \cite{Di2026} obtained a pointwise sparse bound in terms of input percentile oscillations under the conditions $T1\in L^\infty$ and $T^*1\in\mathrm{BMO}$. This bound also implies the mean oscillation sparse estimate above. Our main estimate was obtained independently using a random dyadic representation and a stopping argument for the bilinear pairing. The estimate separates a term involving mean oscillations in both entries from two mixed terms governed by $T1$ and $T^*1$.

\begin{thm}\label{sparse of linear operator}
Let $T$ be a Calder\'on-Zygmund operator, then for all compactly supported functions $f_1,f_2\in L^\infty(\R^d)$, there exists an $\eta$-sparse dyadic family $\mathcal{S}=\mathcal{S}(f_1,f_2)$ such that
\begin{align*}
&\abs{\ave{Tf_1,f_2}}\lesssim\sum_{Q\in\mathcal{S}}\abs{Q}\Big(\Omega(f_1;Q)\Omega(f_2;Q)+\Norm{T1}{\BMO}\ave{\abs{f_1}}_Q\Omega(f_2;Q)\\
&\hspace*{7cm}+\Norm{T^*1}{\BMO}\Omega(f_1;Q)\ave{\abs{f_2}}_Q\Big).
\end{align*}
Moreover, if $T1\in L^\infty$, the term $\ave{\abs{f_1}}_Q$ can be replaced by $\Omega(f_1;Q)$, and $\Norm{T1}{\mathrm{BMO}}$ can be replaced by $\Norm{T1}{L^\infty}$. Similarly, if $T^*1\in L^\infty$, the corresponding replacement applies to the term involving $\ave{\abs{f_2}}_Q$, with $\Norm{T^*1}{\mathrm{BMO}}$ replaced by $\Norm{T^*1}{L^\infty}$.
\end{thm}
When both testing symbols have bounded representatives, Theorem \ref{sparse of linear operator} gives a sparse bound with mean oscillations in both entries. A cutoff test shows that the boundedness of both symbols is also necessary.
\begin{cor}\label{both oscillatory}
   Let $T$ be a Calder\'on-Zygmund operator, then for all compactly supported functions $f_1,f_2\in L^\infty(\R^d)$, there exists an $\eta$-sparse dyadic family $\mathcal{S}$ such that 
\[\abs{\ave{Tf_1,f_2}}\lesssim\sum_{Q\in\mathcal{S}}\abs{Q}\Omega(f_1;Q)\Omega(f_2;Q)\]
if and only if $T1,T^*1\in L^\infty$ simultaneously.
\end{cor}
We next consider oscillation in the input entry alone. The resulting sparse bound characterizes the boundedness of $T1$, which can also be detected by the growth of the unweighted local $L^p$ moments of $Tf$.
\begin{thm}\label{thm2}
   Let $T$ be a Calder\'on-Zygmund operator, then the following statements are equivalent:
   \begin{enumerate}
        \item $T1\in L^\infty$.
      \item For every pair of compactly supported $f_1,f_2\in L^\infty(\R^d)$, there exists an $\eta$-sparse dyadic family $\mathcal{S}$ such that 
      \[\abs{\ave{Tf_1,f_2}}\lesssim\sum_{Q\in\mathcal{S}}\abs{Q}\Omega(f_1;Q)\ave{\abs{f_2}}_Q;\] 
      \item $T$ locally maps $\BMO\to \exp L$, that is, for every cube $R$, when $f\in\BMO$ and $\supp f\subset R$, then 
      \[
      \Norm{Tf}{\exp L(R,\ud x/\abs{R})}\lesssim\Norm{f}{\BMO}.
      \]
     Equivalently, the local $L^p$ moments grow at most linearly in $p$: for all $p\geq 2$,
   \[
   \Big(\frac{1}{\abs{R}}\int_R\abs{Tf(x)}^p\ud x\Big)^{1/p}\lesssim p\Norm{f}{\BMO},
   \]
   where the implicit constant is independent of $R,f,p$.
   \end{enumerate}
\end{thm}
The two formulations in $(3)$ are equivalent by the moment characterization of exponential Orlicz spaces \cite[Section 2.2.4]{T1993}, namely
\[
\Norm{f}{\exp(L^{1/\sigma})(R,\ud x/\abs{R})}
\asymp_\sigma
\sup_{p\ge 2}p^{-\sigma}\Big(\frac{1}{\abs{R}}\int_R\abs{f(x)}^p\ud x\Big)^{1/p}.
\]
The oscillatory sparse bound in $(2)$ also yields the sufficiency of Lerner's pointwise Sobolev characterization \cite[Theorem~1.2]{L2026}
\begin{cor}\cite[Theorem 1.2]{L2026}\label{lerner thm}
   Let $T$ be a Calder\'on-Zygmund operator with $T1\in L^\infty$. Then the Sobolev-type inequality
   \[
   \abs{Tf(x)}\lesssim I_1(\abs{\nabla f})(x)
   \]
   holds for every $f\in C_c^1(\R^n), n\geq 2$, and for almost every $x\in\R^n$.
\end{cor}

The remainder of the paper is organized as follows. In Section 2, we introduce the necessary notation and establish several preliminary results. Section 3 is devoted to the proof of the sparse domination for the model operator, which together with the representation theorem yields Theorem \ref{sparse of linear operator}.  In Section 4, we will prove the characterizations of Theorem \ref{thm2} and recover Lerner's Sobolev-type inequality. 
\section{Preliminaries}
\subsection{Calder\'on-Zygmund operators}
\begin{defn}
A function $K(x,y):\R^d\times\R^d\backslash\{(x,y):x=y\}\to\C$ is called a Calder\'on-Zygmund kernel if it satisfies the size condition:
\begin{align*}
   \abs{K(x,y)}\lesssim\frac{1}{\abs{x-y}^d},
\end{align*}
and there exists some $\alpha\in(0,1]$ such that
\begin{align*}
    \abs{K(x,y)-K(x',y)}+\abs{K(y,x)-K(y,x')}\lesssim\frac{\abs{x-x'}^\alpha}{\abs{x-y}^{d+\alpha}},\text{ whenever }\abs{x-x'}\leq\frac{\abs{x-y}}{2}.
\end{align*}
\end{defn}
\begin{defn}
   We say that $T$ is a Calder\'on-Zygmund operator if $T$ is $L^2(\R^d)$ bounded and there exists a Calder\'on-Zygmund kernel $K$ such that
\[
\ave{Tf_1,f_2}=\int_{\R^d}\int_{\R^d}K(x,y)f_1(y)f_2(x)\ud y\ud x,\,\text{ for }\supp(f_1)\cap\supp(f_2)=\emptyset.
\]  
We also define its adjoint operator $T^*$ by $\ave{Tf_1,f_2}=\ave{T^*f_2,f_1}$. The distributions $T1$ and $T^*1$ are understood modulo constants. The condition $T1\in L^\infty$ means that $T1$ admits a bounded representative.
 \end{defn}
\subsection{Dyadic grids and Haar functions}
Begin with the standard dyadic grid $\mathcal{D}^0$ in $\R^d$, namely
\[
\mathcal{D}^0:=\bigcup_{k\in\Z}\left\{2^{-k}([0,1)^d+m):m\in\Z^d\right\},
\] 
then we define the shifted lattice $\mathcal{D}^\sigma$ by 
\begin{align*}
  \mathcal{D}^\sigma:=\Big\{L+\sigma:=L+\sum_{i:2^{-i}<\ell(L)}2^{-i}\sigma_i: L\in\mathcal{D}^0\Big\}\text{ for }\sigma=(\sigma_i)_{i\in\Z}\in(\{0,1\}^d)^{\Z}.
\end{align*}
Given $I\in\mathcal{D}^\sigma$, we use the following notations through this paper.
\begin{enumerate}
    \item $\ell(I)$ denotes the sidelength of $I$;
    \item $I^{(k)}$ denotes the $k$-th parent of $I$, $I\subset I^{(k)}$ and $\ell(I^{(k)})=2^k\ell(I)$;
    \item $\mathrm{ch}(I)$ is the collection of the children of $I$, i.e., $\mathrm{ch}(I)=\{J\in\mathcal{D}^\sigma:J^{(1)}=I\}$;
    \item $\Exp_{\sigma}$ denotes the expectation with respect to $\sigma\in(\{0,1\}^d)^{\Z}$.
\end{enumerate}
\begin{defn}
    Let $h_I^0,h_I^1$ be standard Haar functions on interval $I\in\R$, namely
    $$h_I^0=\abs{I}^{-1/2}1_I\text{ and }h_I^1=\abs{I}^{-1/2}(1_{I^{-}}-1_{I^{+}}),$$
    where $I^-,I^+$ denote the left and right halves of $I$ respectively. The Haar functions in higher dimension are defined by 
    \[
    h_I^\eta:=h_{I_1}^{\eta_1}\otimes\cdots\otimes h_{I_d}^{\eta_d},\,\text{ for cube }I=\prod_{j=1}^{d}I_j\in\R^d,\eta=(\eta_1,\dots,\eta_d)\in\{0,1\}^d.
    \] 
    Throughout, we suppress cancellative Haar signatures and the corresponding dependence of the coefficients in the notation. Summation over these signatures is understood in the definitions of shifts and paraproducts. In each paraproduct term, the two Haar functions associated with the same cube carry the same signature.
    Sometimes, we need certain more general "Haar functions". Given $I,J\in\mathcal{D}^\sigma$ with $\ell(I)=\ell(J)$, let $H_{I,J}$ be the function satisfying
    \begin{enumerate}
      \item $H_{I,J}$ is supported on $I\cup J$ and constant on their children $\mathrm{ch}(I)\cup\mathrm{ch}(J)$;
      \item $\abs{H_{I,J}}\leq\abs{I}^{-1/2}$;
      \item $\int H_{I,J}\ud x=0$.
    \end{enumerate}
\end{defn}
\subsection{Model operators and dyadic representation theorem}
\begin{defn}\label{def:shift}
   For $k\in\No$, we say that $Q_{k,\sigma}$ is a \textit{dyadic shift of $k$-complexity}, if it has the form 
   \[
   \ave{Q_{k,\sigma}f_1,f_2}=\sum_{\substack{K\in\mathcal{D}^\sigma}}\sum_{\substack{I,J\in\mathcal{D}\\I^{(k)}=J^{(k)}=K}}a_{IJK}\ave{f_1,H_{I,J}}\ave{f_2,h_J},\text{where }\abs{a_{IJK}}\lesssim\frac{\abs{I}^{1/2}\abs{J}^{1/2}}{\abs{K}},
   \]
   or symmetrical form (i.e. the role of $f_1$ is replaced by $f_2$).
\end{defn}
\begin{defn}
    We say that $\pi_{b,\sigma}$ is a \textit{ dyadic paraproduct}, if 
    \[
    \ave{\pi_{b,\sigma}f_1,f_2}=\sum_{I\in\mathcal{D}^\sigma}\ave{b,h_I}\ave{f_1}_I\ave{f_2,h_I},
    \]
    and its adjoint operator is defined by $\ave{\pi_{b,\sigma}^*f_1,f_2}:=\ave{f_1,\pi_{b,\sigma}f_2}$
\end{defn}
We now recall the dyadic representation theorem, which represents the bilinear form $\ave{Tf_1,f_2}$ as a sum of dyadic shifts and dyadic paraproducts above in the expectation sense. Here, we use the modern form of the dyadic representation theorem due to \cite[Theorem 2.26]{AMV2022}.
\begin{thm}\label{representation theorem}
  Let $T$ be a Calder\'on-Zygmund operator, then 
\begin{align*}
  \ave{Tf_1,f_2}=\Exp_{\sigma}\left[\sum_{k=0}^\infty2^{-k\alpha}\ave{Q_{k,\sigma}f_1,f_2}+\ave{\pi_{b_1,\sigma} f_1,f_2}+\ave{\pi_{b_2,\sigma}^* f_1,f_2}\right],
\end{align*}
where $b_1=T1,b_2=T^*1$ and $b_i\in\BMO$ for $i=1,2$, respectively. Here $Q_{k,\sigma}$ denotes a sum of a dimensionally bounded number of shifts of the forms in Definition \ref{def:shift}, with coefficient bounds uniform in $k$ and $\sigma$. 
\end{thm}
\section{Sparse domination of bilinear forms}
The usual sparse domination of bilinear form by local averages had been considered in \cite{LMO2019}. Here, however, the mean-zero property of the Haar function allows us to sharpen the usual sparse bound by replacing size averages with local mean oscillations. Moreover, the representation theorem makes it possible to identify which components contribute only a single oscillation and which ones allow oscillations on both sides. In the stopping conditions below, ratios of the form $0/0$ are interpreted as zero.
\subsection{Oscillatory sparse domination for shifts}
\begin{lem}\label{sparse of shift}
   Let $Q_{k,\sigma}$ be a dyadic shift of $k$-complexity, then for all compactly supported functions $f_1,f_2\in L^\infty(\R^d)$, there exists a sparse family $\mathcal{S}\subset \mathcal{D}^\sigma$ such that 
\[
\abs{\ave{Q_{k,\sigma}f_1,f_2}}\lesssim (k+1)\sum_{Q\in\mathcal{S}}\abs{Q}\Omega(f_1;Q)\Omega(f_2;Q).
\]
\end{lem}
\begin{proof}
   We suppress $\sigma$ throughout the proof since it is fixed. We first consider a finite
truncation of $Q_k$, retaining the same notation. Thus, all but finitely
many coefficients $a_{IJK}$ are set to zero. Choose a sufficiently large
dyadic scale so that every principal cube $K$ is contained in some cube at
this scale. These top cubes are pairwise disjoint, and we argue separately
on each of them. The coefficients of the truncated shift remain fixed
throughout the stopping-time arguments.

Fix one such top cube $Q_0$ and set
\[
\Lambda_{Q_0}(f_1,f_2)
:=
\big\langle
Q_k\big((f_1-\langle f_1\rangle_{Q_0})\mathbf 1_{Q_0}\big),
(f_2-\langle f_2\rangle_{Q_0})\mathbf 1_{Q_0}
\big\rangle.
\]
Since the Haar functions have mean zero and $K\subseteq Q_0$, subtracting these averages does not change the value.
The estimates below are uniform in the truncation and in the choice of
$Q_0$.
 
    Let $\mathcal{E}$ be the collection of all maximal dyadic cubes $Q\in\mathcal{D},Q\subset Q_0$ such that 
    \[
    \max\left(\frac{\ave{\abs{f_1-\ave{f_1}_{Q_0}}}_Q}{\Omega(f_1;Q_0)},\frac{\ave{\abs{f_2-\ave{f_2}_{Q_0}}}_Q}{\Omega(f_2;Q_0)}\right)>C_0,
    \]
    where $C_0$ is so large that 
    \[
    \sum_{Q\in\mathcal{E}}\abs{Q}\leq\frac{1}{2}\abs{Q_0}.
    \]
    Then $\abs{Q_0\backslash\cup_{Q\in\mathcal{E}}Q}\geq1/2\abs{Q_0}$ and we take $Q_0$ to be the first cube in $\mathcal{S}$. Let
    \[
    \mathcal{G}=\left\{Q\in\mathcal{D}:,Q\nsubseteq Q'\text{ for every }Q'\in\mathcal{E}\right\},
    \]
    and $\mathcal{D}(Q)$ denotes $\{K\in\mathcal{D}:K\subset Q\}$. Then we decompose the bilinear form $\Lambda_{Q_0}(f_1,f_2)$ into the following two parts:
    \begin{align*}
        &\uppercase\expandafter{\romannumeral+1}+\uppercase\expandafter{\romannumeral+2}\\
        &:=\left(\sum_{K\in\mathcal{G}}+\sum_{Q\in\mathcal{E}}\sum_{K\in\mathcal{D}(Q)}\right)\sum_{\substack{I,J\in\mathcal{D}\\I^{(k)}=J^{(k)}=K}}a_{IJK}\ave{(f_1-\ave{f_1}_{Q_0})1_{Q_0},H_{I,J}}\ave{(f_2-\ave{f_2}_{Q_0})1_{Q_0},h_J}
    \end{align*}
    Applying the Calder\'on-Zygmund decomposition to $(f_i-\ave{f_i}_{Q_0})1_{Q_0},i=1,2$ with respect to $\mathcal{E}$ gives that for $i=1,2$,
    \begin{align*}
    (f_i-\ave{f_i}_{Q_0})1_{Q_0}=g_i+b_i=g_i+\sum_{Q\in\mathcal{E}}b_{i,Q},\quad b_{i,Q}=\left(f_i-\ave{f_i}_Q\right)1_Q,
     \end{align*}
     where the good and the bad functions satisfy
     \[
     \norm{g_i}{L^\infty}\lesssim\Omega(f_i;Q_0),\,\int b_{i,Q}=0,\,\norm{b_{i,Q}}{L^1}\lesssim\abs{Q}\Omega(f_i;Q_0).
     \]
     Hence it reduces $\uppercase\expandafter{\romannumeral+1}$ to estimate the following four parts:
     \[
    \uppercase\expandafter{\romannumeral+1}(g_1,g_2),\,\uppercase\expandafter{\romannumeral+1}(g_1,b_2),\,\uppercase\expandafter{\romannumeral+1}(b_1,g_2),\,\uppercase\expandafter{\romannumeral+1}(b_1,b_2).
     \]
     The part with two good functions can be directly estimated by the $L^2$ boundedness of the modified shift \cite[Proposition 3.10]{AMV2022},
     \[
     \abs{\uppercase\expandafter{\romannumeral+1}(g_1,g_2)}\lesssim(1+k)^{1/2}\norm{g_1}{L^2}\norm{g_2}{L^2}\lesssim (k+1)\Omega(f_1;Q_0)\Omega(f_2;Q_0)\abs{Q_0}.
     \]
      When it comes to the bad function $b_2$, to guarantee $\ave{b_{2,Q},h_J}\neq 0$ with $J^{(k)}=K\in\mathcal{G}$ and $Q\in\mathcal{E}$, we must have 
      \[
       Q\subsetneq K\quad\text{and}\quad\ell(Q)\geq 2^{-k}\ell(K).
      \] 
      Therefore,
     \begin{align*}
        \uppercase\expandafter{\romannumeral+1}(g_1,b_2)=\sum_{K\in\mathcal{G}}\sum_{\substack{Q\in\mathcal{E},Q\subsetneq K\\\ell(Q)\geq 2^{-k}\ell(K)}}\sum_{\substack{I,J\in\mathcal{D}\\I^{(k)}=J^{(k)}=K}}a_{IJK}\ave{g_1,H_{I,J}}\ave{b_{2,Q},h_J}.
     \end{align*}
     and then 
\begin{align*}
        \abs{\uppercase\expandafter{\romannumeral+1}(g_1,b_2)}\lesssim\sum_{K\in\mathcal{G}}\sum_{\substack{Q\in\mathcal{E},Q\subsetneq K\\\ell(Q)\geq 2^{-k}\ell(K)}}\norm{g_1}{L^\infty}\norm{b_{2,Q}}{L^1}\lesssim (k+1)\Omega(f_1;Q_0)\Omega(f_2;Q_0)\abs{Q_0}.
     \end{align*}
     Similarly, it holds that $\ave{b_{1,Q},H_{I,J}}=0$ for $I^{(k)}=J^{(k)}=K\in\mathcal{G}$ and $Q\in\mathcal{E}$ unless
     \[
       Q\subsetneq K\quad\text{and}\quad\ell(Q)\geq 2^{-k}\ell(K).
      \]
     Then
     \begin{align*}
        &\abs{\uppercase\expandafter{\romannumeral+1}(b_1,g_2)}=\Babs{\sum_{K\in\mathcal{G}}\sum_{\substack{Q\in\mathcal{E},Q\subsetneq K\\\ell(Q)\geq 2^{-k}\ell(K)}}\sum_{\substack{I,J\in\mathcal{D}\\I^{(k)}=J^{(k)}=K}}a_{IJK}\ave{b_{1,Q},H_{I,J}}\ave{g_{2},h_J}}\\
        &\hspace*{1cm}\lesssim\sum_{K\in\mathcal{G}}\sum_{\substack{Q\in\mathcal{E},Q\subsetneq K\\\ell(Q)\geq 2^{-k}\ell(K)}}\norm{g_2}{L^\infty}\norm{b_{1,Q}}{L^1}\lesssim (k+1)\Omega(f_1;Q_0)\Omega(f_2;Q_0)\abs{Q_0}.
     \end{align*}
     Write the term with two bad parts as 
     \[
     \uppercase\expandafter{\romannumeral+1}(b_1,b_2)=\sum_{K\in\mathcal{G}}\sum_{\substack{Q\in\mathcal{E},Q\subsetneq K\\\ell(Q)\geq 2^{-k}\ell(K)}}\sum_{\substack{Q'\in\mathcal{E},Q'\subsetneq K\\\ell(Q)'\geq 2^{-k}\ell(K)}}\sum_{\substack{I,J\in\mathcal{D}\\I^{(k)}=J^{(k)}=K}}a_{IJK}\ave{b_{1,Q},H_{I,J}}\ave{b_{2,Q'},h_J}.
     \]
     Split $\uppercase\expandafter{\romannumeral+1}(b_1,b_2)=\uppercase\expandafter{\romannumeral+1}^1(b_1,b_2)+\uppercase\expandafter{\romannumeral+1}^2(b_1,b_2)$, where $\uppercase\expandafter{\romannumeral+1}^1, \uppercase\expandafter{\romannumeral+1}^2$ correspond to the following decompositions of the sum, respectively.
     \[
     \sum_{\substack{Q\in\mathcal{E},Q\subsetneq K\\\ell(Q)\geq 2^{-k}\ell(K)}}\sum_{\substack{Q'\in\mathcal{E},Q'\subsetneq K\\\ell(Q)'\geq 2^{-k}\ell(K)}}=\sum_{Q\neq Q'}+\sum_{Q= Q'}.
     \]
     When $Q\neq Q'$, for each nonzero summand, the cancellation of $b_{2,Q'}$ implies $J\subseteq Q'$. Then the part of $H_{I,J}$ supported on $J$ does not contribute, namely 
     \[
     \ave{b_{1,Q},H_{I,J}}=\ave{b_{1,Q},1_IH_{I,J}}+\ave{b_{1,Q},1_{I\neq J}1_JH_{I,J}}=\ave{b_{1,Q},1_IH_{I,J}}.
     \]
     Therefore,
     \begin{align*}
        \abs{\uppercase\expandafter{\romannumeral+1}^1(b_1,b_2)}&\lesssim\sum_{K\in\mathcal{G}}\sum_{\substack{Q\in\mathcal{E},Q\subsetneq K\\\ell(Q)\geq 2^{-k}\ell(K)}}\sum_{\substack{Q'\in\mathcal{E},Q'\subsetneq K\\\ell(Q)'\geq 2^{-k}\ell(K)}}\frac{1}{\abs{K}}\norm{b_{1,Q}}{L^1}\norm{b_{2,Q'}}{L^1}\\
&\lesssim (k+1)\Omega(f_1;Q_0)\Omega(f_2;Q_0)\abs{Q_0}.
     \end{align*}
     On the other hand,
     \begin{align*}
    \uppercase\expandafter{\romannumeral+1}^2(b_1,b_2)&=\sum_{K\in\mathcal{G}}\sum_{\substack{Q\in\mathcal{E},Q\subsetneq K\\\ell(Q)\geq 2^{-k}\ell(K)}}\sum_{\substack{I,J\in\mathcal{D}\\I^{(k)}=J^{(k)}=K}}a_{IJK}\ave{b_{1,Q},H_{I,J}}\ave{b_{2,Q},h_J}\\
&=\sum_{Q\in\mathcal{E}}\sum_{\substack{K\in\mathcal{G},Q\subsetneq K\\\ell(Q)\geq 2^{-k}\ell(K)}}\sum_{\substack{I,J\in\mathcal{D}\\I^{(k)}=J^{(k)}=K}}a_{IJK}\ave{b_{1,Q},H_{I,J}}\ave{b_{2,Q},h_J}.
     \end{align*}
     Note that  
    \[
    \uppercase\expandafter{\romannumeral+2}+\uppercase\expandafter{\romannumeral+1}^2(b_1,b_2)=\sum_{Q\in\mathcal{E}}\Lambda_Q(f_1,f_2),
    \]
    iterating this estimate produces a $1/2$-sparse family below $Q_0$.
To justify the iteration, the $L^2$ boundedness of the shift gives that for each cube $Q$,
\[
\begin{aligned}
|\Lambda_Q(f_1,f_2)|
&\lesssim
\sqrt{k+1}\,
\|(f_1-\langle f_1\rangle_Q)\mathbf 1_Q\|_{L^2}
\|(f_2-\langle f_2\rangle_Q)\mathbf 1_Q\|_{L^2}\\
&\lesssim
\sqrt{k+1}\,
\|f_1\|_{L^\infty}\|f_2\|_{L^\infty}|Q|.
\end{aligned}
\]
The total measure of the stopping cubes in generation $m$ is at most
$2^{-m}\abs{Q_0}$, which tends to zero as $m\to\infty$.

Taking the union of the sparse families obtained under the disjoint top
cubes preserves $1/2$-sparseness. Consequently, every finite truncation
satisfies
\[
|\langle Q_k f_1,f_2\rangle|
\lesssim
(k+1)
\sup_{\substack{\mathcal S\subseteq\mathcal D\\
                 \mathcal S\ \text{is }1/2\text{-sparse}}}
\sum_{Q\in\mathcal S}|Q|\Omega(f_1;Q)\Omega(f_2;Q),
\]
with a constant independent of the truncation.
The supremum is finite, since
$\Omega(f_i;Q)\le 2\langle |f_i|\rangle_Q$, and sparseness together with
the $L^2$ boundedness of the Hardy--Littlewood maximal operator bounds
each such sum by a constant multiple of
$\|f_1\|_{L^2}\|f_2\|_{L^2}$.
Passing to the limit in the pairings by $L^2$ convergence yields the
same estimate for the original shift. If the supremum is positive,
choose a sparse family whose sum exceeds half of it. If it is zero,
the pairing vanishes. This proves the lemma.
\end{proof}
\subsection{Oscillatory sparse domination for paraproducts}
\begin{lem}\label{sparse of paraproduct}
   Let $\pi_{b,\sigma}$ be a dyadic paraproduct with $b\in\BMO$, then for all compactly supported functions $f_1,f_2\in L^\infty(\R^d)$, there exists a sparse family $\mathcal{S}\subset\mathcal{D}^\sigma$ such that 
\[
\abs{\ave{\pi_{b,\sigma}f_1,f_2}}\lesssim \Norm{b}{\BMO}\sum_{Q\in\mathcal{S}}\abs{Q}\ave{\abs{f_1}}_Q\Omega(f_2;Q).
\]
If, in addition, $b\in L^\infty$, then we have 
\[
\abs{\ave{\pi_{b,\sigma}f_1,f_2}}\lesssim \Norm{b}{L^\infty}\sum_{Q\in\mathcal{S}}\abs{Q}\Omega(f_1;Q)\Omega(f_2;Q).
\]
\end{lem}
\begin{proof}
    We first restrict the paraproduct to cubes whose side lengths are at most a fixed dyadic scale, and work separately on each top cube at that scale, namely
   \[
   \sum_{\substack{I\in\mathcal{D}\\I\subset Q_0}}\ave{b,h_I}\ave{f_1}_I\ave{f_2,h_I}.
   \]
    Let $\mathcal{E}$ be the collection of all maximal dyadic cubes $Q\in\mathcal{D},Q\subset Q_0$ such that 
    \[
    \max\left(\frac{\ave{\abs{f_1}}_Q}{\ave{\abs{f_1}}_{Q_0}},\frac{\ave{\abs{f_2-\ave{f_2}_{Q_0}}}_Q}{\Omega(f_2;Q_0)}\right)>C_0,
    \]
    where we choose $C_0$ is sufficiently large such that $\sum_{Q\in\mathcal{E}}\abs{Q}\leq1/2\abs{Q_0}$.
    As the proof of Lemma \ref{sparse of shift}, we define $\mathcal{G}=\{Q\in\mathcal{D}:Q\subset Q_0,Q\nsubseteq Q'\text{ for every }Q'\in\mathcal{E}\}$.
    Then, by the same reasoning, it suffices to prove
    \[
    \left\vert\sum_{\substack{I\in\mathcal{G}}}\ave{b,h_I}\ave{f_1}_I\ave{f_2,h_I}\right\vert\lesssim\Norm{b}{\BMO}\ave{\abs{f_1}}_{Q_0}\Omega(f_2;Q_0)\abs{Q_0}.
    \]
The cancellation of the Haar function implies that the expression inside the absolute value on the left can be rewritten as
    \[
    \Lambda(f_1,f_2):=\sum_{\substack{I\in\mathcal{G}}}\ave{b,h_I}\ave{f_1}_I\ave{f_2-\ave{f_2}_{Q_0},h_I}.
    \]
    Let 
    \[
    F_1(x)=
    \begin{cases}
      \ave{f_1}_Q &\text{ if }x\in Q\in\mathcal{E},\\
      f_1(x) &\text{ if }x\in Q_0\backslash\cup_{Q\in\mathcal{E}}Q,
    \end{cases}
    ,\,
    F_2(x)=
    \begin{cases}
      \ave{f_2}_Q-\ave{f_2}_{Q_0} &\text{ if }x\in Q\in\mathcal{E},\\
      f_2(x)-\ave{f_2}_{Q_0} &\text{ if }x\in Q_0\backslash\cup_{Q\in\mathcal{E}}Q,
    \end{cases}
    \]
    and both functions are extended by zero outside $Q_0$. Then we have 
    \[
    \Norm{F_1}{L^\infty}\lesssim\ave{\abs{f_1}}_{Q_0}\quad\text{ and }\quad\Norm{F_2}{L^\infty}\lesssim\Omega(f_2;Q_0).
    \]
    For any $I\in\mathcal{G}$, if $I\cap Q\neq \emptyset$ for some $Q\in\mathcal{E}$, then it must be $Q\subsetneq I$. Hence for any $I\in\mathcal{G}$, it holds that 
    \[
    \ave{f_1}_I=\ave{F_1}_I\,\text{ and }\,\ave{ f_2-\ave{f_2}_{Q_0},h_I}=\ave{ F_2,h_I}.
    \]
    Therefore, by the Carleson condition of $b\in\BMO$,
     \begin{align*}
     \abs{\Lambda(f_1,f_2)}&=\Babs{\sum_{\substack{I\in\mathcal{G}}}\ave{b,h_I}\ave{F_1}_I\ave{F_2,h_I}}\leq\Norm{F_2}{L^2}\BNorm{\sum_{\substack{I\in\mathcal{G}}}\ave{b,h_I}\ave{F_1}_Ih_I}{L^2}\\
     &\lesssim\Norm{F_2}{L^2}\Big(\sum_{I\subset Q_0}\abs{\ave{b,h_I}}^2\abs{\ave{F_1}_I}^2\Big)^{1/2}\lesssim\Norm{b}{\BMO}\norm{F_1}{L^2}\norm{F_2}{L^2}\\
     &\lesssim\Norm{b}{\BMO}\abs{Q_0}\norm{F_1}{L^\infty}\norm{F_2}{L^\infty}\lesssim\Norm{b}{\BMO}\abs{Q_0}\ave{\abs{f_1}}_{Q_0}\Omega(f_2;Q_0).
     \end{align*}

     If, in addition, we have $b\in L^\infty$, then write
     \begin{align*}
   &\sum_{\substack{I\in\mathcal{D}\\I\subset Q_0}}\ave{b,h_I}\ave{f_1}_I\ave{f_2,h_I}=\Lambda_{Q_0}(f_1,f_2)+\ave{f_1}_{Q_0}E(Q_0)\\
   &:=\sum_{\substack{I\in\mathcal{D}\\I\subset Q_0}}\ave{b,h_I}\ave{f_1-\ave{f_1}_{Q_0}}_I\ave{f_2,h_I}+\ave{f_1}_{Q_0}\int_{Q_0}(b(x)-\ave{b}_{Q_0})(f_2(x)-\ave{f_2}_{Q_0})\ud x.
   \end{align*}
   For $\Lambda_{Q_0}(f_1,f_2)$, The argument is the same as in the case $b\in \mathrm{BMO}$, except that $f_1$ is replaced by $f_1-\ave{f_1}_{Q_0}$. Thus, let $\mathcal{E}$ be the collection of all maximal dyadic cubes $Q\in\mathcal{D},Q\subset Q_0$ such that 
    \[
    \max\left(\frac{\ave{\abs{f_1-\ave{f_1}_{Q_0}}}_Q}{\Omega(f_1;Q_0)},\frac{\ave{\abs{f_2-\ave{f_2}_{Q_0}}}_Q}{\Omega(f_2;Q_0)}\right)>C_0,
    \]
    and
     $\mathcal{G}=\{Q\in\mathcal{D}:Q\subset Q_0,Q\nsubseteq Q'\text{ for every }Q'\in\mathcal{E}\}$,
    we have 
    \[
    \left\vert\sum_{\substack{I\in\mathcal{G}}}\ave{b,h_I}\ave{f_1-\ave{f_1}_{Q_0}}_I\ave{f_2,h_I}\right\vert\lesssim\Norm{b}{\BMO}\Omega(f_1;Q_0)\Omega(f_2;Q_0)\abs{Q_0}.
    \]
    We write the remaining term as 
\begin{align*}
&\sum_{Q\in\mathcal{E}}\sum_{I\in\mathcal{D}(Q)}\ave{b,h_I}\ave{f_1-\ave{f_1}_{Q_0}}_I\ave{f_2,h_I}\\
&\hspace*{1cm}=\sum_{Q\in\mathcal{E}}\sum_{I\in\mathcal{D}(Q)}\ave{b,h_I}\ave{f_1-\ave{f_1}_{Q}}_I\ave{f_2,h_I}+\sum_{Q\in\mathcal{E}}(\ave{f_1}_Q-\ave{f_1}_{Q_0})E(Q)\\
&\hspace*{1cm}=\sum_{Q\in\mathcal{E}}\Lambda_Q(f_1,f_2)+\sum_{Q\in\mathcal{E}}(\ave{f_1}_Q-\ave{f_1}_{Q_0})E(Q).
\end{align*}
Next, since $b\in L^\infty$, we have 
\begin{align*}
\Babs{\sum_{Q\in\mathcal{E}}(\ave{f_1}_Q-\ave{f_1}_{Q_0})E(Q)}&\lesssim\Norm{b}{L^\infty}\Omega(f_1;Q_0)\sum_{Q\in\mathcal{E}}\int_Q\abs{f_2-\ave{f_2}_Q}\\
&\lesssim\Norm{b}{L^\infty}\Omega(f_1;Q_0)\Omega(f_2;Q_0)\abs{Q_0}.
\end{align*}
Iterating this estimate gives the required sparse bound for $\Lambda_{Q_0}$. The recursive remainder tends to zero by the $L^2$ boundedness of the paraproduct and the boundedness of $f_1,f_2$. Regarding the average remaining from the first step $\ave{f_1}_{Q_0}E(Q_0)$, note that 
\[
\abs{\ave{f_1}_{Q_0}E(Q_0)}\lesssim\frac{\Norm{b}{L^\infty}}{\abs{Q_0}}\Norm{f_1}{L^1(Q_0)}\Norm{f_2}{L^1(Q_0)},
\]
which tends to zero as we take increasingly large top cubes $Q_0$. Summing over the disjoint top cubes preserves sparseness. Moreover,
\[
\sum_{Q_0}\babs{\langle f_1\rangle_{Q_0}E(Q_0)}\to0
\]
as their common scale tends to infinity. Passing to the limit in $L^2$ and applying the supremum argument from Lemma \ref{sparse of shift} completes the proof.
 \end{proof}
 \begin{rem}
   The same result holds for $\pi_b^*$, except that $f_1$ is replaced by $f_2$. 
 \end{rem}
 \subsection{Proof of Theorem \ref{sparse of linear operator}}
   From Theorem \ref{representation theorem}, Lemma \ref{sparse of shift}, Lemma \ref{sparse of paraproduct} and 
   \[
   \sum_{k\geq 0}(k+1)2^{-k\alpha}<\infty,
   \]
   we obtain
   \begin{align*}
   &\abs{\ave{Tf_1,f_2}}\lesssim\sup_{\substack{\sigma\in(\{0,1\}^d)^{\Z}\\\mathcal{S}\subset \mathcal{D}^\sigma\,\eta\text{-sparse}}}\sum_{Q\in\mathcal{S}}\abs{Q}\Big(\Omega(f_1;Q)\Omega(f_2;Q)\\
   &\hspace*{5cm}+\Norm{T1}{\BMO}\ave{\abs{f_1}}_Q\Omega(f_2;Q)+\Norm{T^*1}{\BMO}\Omega(f_1;Q)\ave{\abs{f_2}}_Q\Big).
    \end{align*}
    By the boundedness of maximal operators and sparseness, it holds that
    \begin{align*}
    &\sup_{\substack{\sigma\in(\{0,1\}^d)^{\Z}\\\mathcal{S}\subset \mathcal{D}^\sigma\,\eta\text{-sparse}}}\sum_{Q\in\mathcal{S}}\abs{Q}\Big(\Omega(f_1;Q)\Omega(f_2;Q)\\
   &\hspace*{3cm}+\Norm{T1}{\BMO}\ave{\abs{f_1}}_Q\Omega(f_2;Q)+\Norm{T^*1}{\BMO}\Omega(f_1;Q)\ave{\abs{f_2}}_Q\Big)=C<\infty.
   \end{align*}
   If $C=0$, the desired result is immediate. Otherwise one can select some $\sigma_0$ and $\eta$-sparse family $\mathcal{S}$ such that 
   \begin{align*}
   &\sum_{Q\in\mathcal{S}\subset \mathcal{D}^{\sigma_0}}\abs{Q}\Big(\Omega(f_1;Q)\Omega(f_2;Q)\\
   &\hspace*{2cm}+\Norm{T1}{\BMO}\ave{\abs{f_1}}_Q\Omega(f_2;Q)+\Norm{T^*1}{\BMO}\Omega(f_1;Q)\ave{\abs{f_2}}_Q\Big)>\frac{C}{2},
   \end{align*}
   hence we obtain 
\begin{align*}
&\abs{\ave{Tf_1,f_2}}\lesssim\sum_{Q\in\mathcal{S}\subset \mathcal{D}^{\sigma_0}}\abs{Q}\Big(\Omega(f_1;Q)\Omega(f_2;Q)\\
   &\hspace*{3cm}+\Norm{T1}{\BMO}\ave{\abs{f_1}}_Q\Omega(f_2;Q)+\Norm{T^*1}{\BMO}\Omega(f_1;Q)\ave{\abs{f_2}}_Q\Big).
\end{align*}
If $T1\in L^\infty$ or $T^*1\in L^\infty$, we apply the corresponding estimate in Lemma \ref{sparse of paraproduct} and repeat the same supremum argument. The resulting double-oscillation terms can then be absorbed into the first term.

 \section{Bounded testing conditions and application}
 \subsection{Converses to the oscillatory sparse bounds}
 In this subsection, we give the proof of $(2)\Rightarrow(1)$ in Theorem \ref{thm2}. The proof of $(1)\Rightarrow(2)$  has been already covered by Theorem \ref{sparse of linear operator}, so it suffices to verify the following lemma, whose proof follows Lerner's cutoff argument for the necessity direction, obtaining the uniform bound directly from the sparse form.
 \begin{lem}\label{bound-sparse}
   Let $T$ be a Calder\'on-Zygmund operator, and for all compactly supported functions $f_1,f_2\in L^\infty(\R^d)$, there exists a sparse family $\mathcal{S}$ such that 
\[
\abs{\ave{Tf_1,f_2}}\lesssim\sum_{Q\in\mathcal{S}}\abs{Q}\Omega(f_1;Q)\ave{\abs{f_2}}_Q.
\]
Then $T1\in L^\infty$.
 \end{lem}
 \begin{proof}
    Let $\theta\in C_c^\infty$ with $\theta(x)\equiv 1$ for $\abs{x}\leq 1$, and $\theta_R(x):=\theta(x/R)$. For any cube $Q$,
   \begin{align*}
   \omega(\theta_R;Q)=\frac{1}{\abs{Q}}\int_Q\abs{\theta_R(x)-\ave{\theta_R}_Q}\ud x&\lesssim\frac{1}{\abs{Q}^2}\int_{Q\times Q}\abs{\theta_R(x)-\theta_R(y)}\ud x\ud y\\
&\lesssim\min\Big\{\frac{\ell(Q)}{R},\frac{R^d}{\abs{Q}}\Big\},
   \end{align*}
   where the first bound is from $\Norm{\theta_R}{L^1}\lesssim R^d$ and another one is by $\Norm{\nabla \theta_R}{L^\infty}\lesssim R^{-1}$. Test the oscillatory sparse domination with $f_1=\theta_R$ and $f_2=\varphi\in \mathcal{D}_0(\R^d)$, where $\mathcal{D}_0(\R^d)$ denotes the space of all smooth compactly supported functions on $\R^d$ with mean zero, then 
   \[
   \abs{\ave{T\theta_R,\varphi}}\lesssim\sum_{Q\in\mathcal{S}}\abs{Q}\min\Big\{\frac{\ell(Q)}{R},\frac{R^d}{\abs{Q}}\Big\}\ave{\abs{\varphi}}_Q.
   \]
   Split the cubes in $\mathcal{S}$ as $\sum_{Q\in\mathcal{S}}=\sum_{k\in \Z}\sum_{\substack{Q\in\mathcal{S}\\2^{k}R<\ell(Q)\leq 2^{k+1}R}}$, then we get that 
   \[
   \abs{\ave{T\theta_R,\varphi}}\lesssim\sum_{k\in \Z}\min\{2^k,2^{-kd}\}\sum_{\substack{Q\in\mathcal{S}\\2^{k}R<\ell(Q)\leq 2^{k+1}R}}\abs{Q}\ave{\abs{\varphi}}_Q.
   \]
   Fix $k$ and $x$, and consider any finite collection of $N$ cubes $Q\in\mathcal S$ satisfying $2^kR<\ell(Q)\le 2^{k+1}R$ and $x\in Q$. Then these $N$ cubes are contained in a larger cube $P$, which is centered at $x$ and $\ell(P)=2^{k+2}R$. The sparseness then gives disjoint subsets $E_Q\subset Q$ with
   \[
   \abs{E_Q}\geq \eta\abs{Q}>\eta 2^{kd}R^d.
   \]
   Therefore, we have 
   \[
   N\eta 2^{kd}R^d\leq \Babs{\bigcup_Q E(Q)}\leq 2^{(k+2)d}R^d,
   \]
   then $N\leq \eta^{-1}4^d$, which means that these cubes have uniformly bounded overlap. Then we write
   \[
   \sum_{\substack{Q\in\mathcal{S}\\2^{k}R<\ell(Q)\leq 2^{k+1}R}}\abs{Q}\ave{\abs{\varphi}}_Q\lesssim\int \abs{\varphi(x)}\sum_{\substack{Q\in\mathcal{S}\\2^{k}R<\ell(Q)\leq 2^{k+1}R}}1_Q(x)\ud x\lesssim\Norm{\varphi}{L^1}.
   \]
   Substituting this estimate into the above, we obtain
   \[
   \abs{\ave{T\theta_R,\varphi}}\lesssim\sum_{k\in \Z}\min\{2^k,2^{-kd}\}\Norm{\varphi}{L^1}\lesssim\Norm{\varphi}{L^1}\,\text{ uniformly in }R.
   \]
   Thus, by \cite[Proposition 2.1]{L2026}, there holds that
   \[
   \abs{\ave{T(1),\phi}}=\lim_{R\to\infty}\abs{\ave{T\theta_R,\phi}}\lesssim\Norm{\varphi}{L^1}.
   \]
   By the Hahn-Banach theorem and the duality $(L^1)^*=L^\infty$, there exists $b\in L^\infty$ such that 
   \[
   \ave{T(1),\varphi}=\int b(x)\varphi(x)\ud x\quad \text{ for all } \varphi\in\mathcal{D}_0(\R^d).
   \]
   Therefore $T1\in L^\infty$.
 \end{proof}
 \begin{rem}
 Similar arguments can imply 
 \[
 \abs{\ave{Tf_1,f_2}}\lesssim\sum_{Q\in\mathcal{S}}\abs{Q}\ave{\abs{f_1}}_Q\Omega(f_2;Q)\,\Longrightarrow \, T^*1\in L^\infty.
 \]
 Moreover, if we have 
 \[
 \abs{\ave{Tf_1,f_2}}\lesssim\sum_{Q\in\mathcal{S}}\abs{Q}\Omega(f_1;Q)\Omega(f_2;Q),
 \]
 then by $\Omega(f;Q)\lesssim \ave{\abs{f}}_Q$, there holds that $T1,T^*1\in L^\infty$ simultaneously.
 \end{rem}
 \subsection{Local moment characterizations}
 Before giving our local moment characterizations, we first state the following simpler result.
 \begin{prop}
  Let $T$ be a Calder\'on-Zygmund operator, then for each cube $R$, every function $f\in\BMO$ with $\supp(f)\subset R$ and all $p\geq 2$, it holds that 
   \begin{align*}
   \Big(\frac{1}{\abs{R}}\int_R\abs{Tf(x)}^p\ud x\Big)^{1/p}\lesssim p^2\Norm{f}{\BMO},
   \end{align*}
   where the implicit constant is independent of $R,f,p$.
 \end{prop}
\begin{proof}
   By the classical $L^p$-norm constant of Calder\'on-Zygmund operators, e.g., see \cite[Theorem 1.3]{P2009}, it holds that 
   \[
   \abs{R}^{-1/p}\Norm{Tf}{L^p(R)}\lesssim p \abs{R}^{-1/p}\Norm{f}{L^p}.
   \]
   Since $\supp(f)\subset R$, then 
\[
\abs{\ave{f}_{3R}}=\frac{1}{\abs{3R\backslash R}}\int_{3R\backslash R}\abs{f-\ave{f}_{3R}}\lesssim\Norm{f}{\BMO}.
\]
Therefore, the John-Nirenberg inequality implies
\[
\Norm{f}{L^p}\lesssim \Big(\int_R\abs{f-\ave{f}_{3R}}^p\Big)^{1/p}+\abs{\ave{f}_{3R}}\abs{R}^{1/p}\lesssim p \Norm{f}{\BMO}\abs{R}^{1/p}.
\]
Then we are done. 
\end{proof}
The main result of this section is that, under the stronger assumption $T1\in L^\infty$, the exponent $p^2$ in the original estimate can be reduced to the sharp exponent $p$. Conversely, we show that this improved estimate also implies $T1\in L^\infty$.
 \begin{thm}
   Let $T$ be a Calder\'on-Zygmund operator and $T1\in L^\infty$, then for each cube $R$, every function $f\in\BMO$ with $\supp(f)\subset R$ and all $p\geq 2$, it holds that 
   \begin{align*}
   \Big(\frac{1}{\abs{R}}\int_R\abs{Tf(x)}^p\ud x\Big)^{1/p}\lesssim p\Norm{f}{\BMO},
   \end{align*}
   where the implicit constant is independent of $R,f,p$. Moreover, this result is sharp, in the sense that the power of $p$ in the estimate cannot be lowered.
 \end{thm}
 \begin{proof}
   If $T1\in L^\infty$, we first assume that $f,g\in L^\infty$ are supported in $R$, then 
   \[
   \abs{\ave{Tf,g}}\lesssim\sum_{Q\in\mathcal{S}}\abs{Q}\Omega(f;Q)\ave{\abs{g}}_Q.
   \]
   Next, we handle the cubes in $\mathcal{S}$ with $\ell(Q)\leq \ell(R)$ and $\ell(Q)> \ell(R)$ respectively.

   \textbf{Case 1,} $\ell(Q)\leq \ell(R)$. Then if $Q\cap R\neq \emptyset$, we have $Q\subset 3R$. Hence 
   \begin{align*}
      \sum_{\substack{Q\in\mathcal{S}\\\ell(Q)\leq \ell(R)}}\abs{Q}\Omega(f;Q)\ave{\abs{g}}_Q&\lesssim\Norm{f}{\BMO}\sum_{\substack{Q\in\mathcal{S}\\\ell(Q)\leq \ell(R)}}\abs{E(Q)}\ave{\abs{g}}_Q\\
      &\lesssim\Norm{f}{\BMO}\int_{3R} Mg\lesssim p\Norm{f}{\BMO}\abs{R}^{1/p}\Norm{g}{L^{p'}},
   \end{align*}
where in the last step we use the fact that $\Norm{M}{L^{p'}\to L^{p'}}\lesssim p$.

\textbf{Case 2,} $\ell(Q)> \ell(R)$. Split the cubes as follows
\[
\sum_{\substack{Q\in\mathcal{S}\\\ell(Q)> \ell(R)}}=\sum_{k\geq 0}\sum_{\substack{Q\in\mathcal{S}\\2^k\ell(R)<\ell(Q)\leq 2^{k+1}\ell(R)}}.
\]
For each fixed $k$, if $Q\cap R\neq\emptyset$, then $Q\subset 2^{k+3}R$. Assuming that there are $N$ such cubes, then
\[
N\eta 2^{kd}\abs{R}\leq \Babs{\cup E(Q)}\leq 2^{(k+3)d}\abs{R}.
\]
Therefore we have $N\lesssim1$ and 
\begin{align*}
    \sum_{\substack{Q\in\mathcal{S}\\\ell(Q)> \ell(R)}}\abs{Q}\Omega(f;Q)\ave{\abs{g}}_Q&\lesssim \sum_{k\geq 0}\sum_{\substack{Q\in\mathcal{S}\\2^k\ell(R)<\ell(Q)\leq 2^{k+1}\ell(R)}}\frac{\Norm{f}{L^1}\Norm{g}{L^1}}{\abs{Q}}\lesssim\frac{\Norm{f}{L^1}\Norm{g}{L^1}}{\abs{R}}.
\end{align*}
Since $\supp(f)\subset R$, then 
\[
\abs{\ave{f}_{3R}}=\frac{1}{\abs{3R\backslash R}}\int_{3R\backslash R}\abs{f-\ave{f}_{3R}}\lesssim\Norm{f}{\BMO}.
\]
As a result, we obtain
\[
\Norm{f}{L^1}\lesssim \int_R\abs{f-\ave{f}_{3R}}+\abs{\ave{f}_{3R}}\abs{R}\lesssim\Norm{f}{\BMO}\abs{R}.
\]
Together with $\Norm{g}{L^1}\leq\abs{R}^{1/p}\Norm{g}{L^{p'}}$, there holds that 
\[
\frac{\Norm{f}{L^1}\Norm{g}{L^1}}{\abs{R}}\lesssim \Norm{f}{\BMO}\abs{R}^{1/p}\Norm{g}{L^{p'}}.
\]
Combining these two cases and applying $L^p$ duality give the desired estimate for $f\in L^\infty$. For general $f\in\mathrm{BMO}$ supported in $R$, approximate $f$ by bounded truncations supported in $R$, with uniformly controlled $\BMO$ seminorms and convergence in $L^2$. The $L^2$ boundedness of $T$ and Fatou's lemma then yield the estimate for $f$.

 As for the sharpness, let $T=H$ be the Hilbert transform, $R=(0,1)$ and $f=1_{R}$, then $\Norm{f}{\BMO}\lesssim1$ and for $x\in(0,1)$,
 \[
 Hf(x)=\frac{1}{\pi}\log\frac{x}{1-x}.
 \]
 Then on $0<x<e^{-p}/2$, we have $\abs{Hf(x)}\geq p/\pi$. Hence
 \[
 \Big(\int_0^1\abs{Hf}^p\Big)^{1/p}\geq \frac{p}{\pi}\Big(\frac{e^{-p}}{2}\Big)^{1/p}\gtrsim p.
 \]
\end{proof}
\begin{rem}
This theorem also shows that the target space in the local $\mathrm{BMO}\to\exp L$ estimate cannot be improved to $\exp(L^{1+\varepsilon})$ for any $\varepsilon>0$.
\end{rem}
\begin{thm}
   Let $T$ be a Calder\'on-Zygmund operator and  for each cube $R$, every function $f\in\BMO$ with $\supp(f)\subset R$ and all $p\geq 2$, we have
   \begin{align}\label{p-norm}
   \Big(\frac{1}{\abs{R}}\int_R\abs{Tf(x)}^p\ud x\Big)^{1/p}\lesssim p\Norm{f}{\BMO},
   \end{align}
   where the implicit constant is independent of $R,f,p$, then there holds that $T1\in L^\infty$. 
 \end{thm}
 \begin{proof}
Suppose that the estimate (\ref{p-norm}) holds, let 
\[
F_N(x)=\sum_{j=1}^N\theta(2^{-j}x),
\]
where $\theta$ is the same as the $\theta$ in the proof of Lemma \ref{bound-sparse}, we claim that $F_N\in\BMO$ and $\Norm{F_N}{\BMO}\lesssim1$. Set
\[ 
H_N(x)=\int_0^N1_{\{\abs{x}\leq 2^t\}}\ud t=\min\{N,(N-\log_2\abs{x})_+\},
\]
then we have 
\[
\abs{F_N(x)-H_N(x)}\leq \sum_{j=1}^N\abs{\theta(2^{-j}x)-1_{\{\abs{x}\leq 2^j\}}}+\Babs{\sum_{j=1}^N\Big(1_{\{\abs{x}\leq 2^j\}}-\int_{j-1}^j1_{\{\abs{x}\leq 2^t\}}\ud t\Big)}\lesssim1.
\]
Hence the claim reduces to $\Norm{H_N}{\BMO}\lesssim1$. For every cube $Q$, since $\Norm{\log_2\abs{\cdot}}{\BMO}\lesssim1$, then there holds that
\begin{align*}
   \int_Q\abs{H_N-\ave{H_N}_Q}&\lesssim\int_Q\abs{\min\{N,(N-\log_2\abs{x})_+\}-\min\{N,(N-\ave{\log_2\abs{x}}_Q)_+\}}\ud x\\
   &\lesssim\int_Q\abs{\log_2\abs{x}-\ave{\log_2\abs{x}}_Q}\ud x\lesssim\abs{Q}.
\end{align*}
So we complete the proof of the claim. Note that the support of $F_N$ is contained in some cube $R_N$ with $\abs{R_N}\sim 2^{Nd}$, applying (\ref{p-norm}) with $f=F_N/N,p=N$ and $R=R_N$ gives 
\[
\BNorm{T(\frac{F_N}{N})}{L^N(R_N)}\lesssim \abs{R_N}^{1/N}\Norm{F_N}{\BMO}\lesssim 1.
\]
Take $\varphi\in \mathcal{D}_0$, when $N$ is sufficiently large, $\supp(\varphi)\subset R_N$. By duality, we have 
\[
\abs{\ave{T(F_N/N),\varphi}}\lesssim\Norm{\varphi}{L^{N/(N-1)}}.
\]
Note that $\ave{T(F_N/N),\varphi}$ is the Ces\`aro average of $\ave{T(\theta(2^{-j}\cdot)),\varphi}$, so it converges to the same limit. Let $N\to\infty$, we have 
\[
\ave{T(F_N/N),\varphi}\to\ave{T1,\varphi}.
\] 
On the other hand, for every $1<p\leq 2$, we have $\abs{\varphi(x)}^p\leq\abs{\varphi(x)}+\abs{\varphi(x)}^2\in L^1$. Thus by dominated convergence theorem,
\[
\Norm{\varphi}{L^p}^p=\int \abs{\varphi}^p\to\int \abs{\varphi}=\Norm{\varphi}{L^1}\quad \text{ as }p\to 1,
\]
and then $\Norm{\varphi}{L^p}=\exp(1/p\log\Norm{\varphi}{L^p}^p)\to \exp(\log\Norm{\varphi}{L^1})=\Norm{\varphi}{L^1}$ as $p\to 1$. Therefore,
\[
\abs{\ave{T1,\varphi}}\lesssim\Norm{\varphi}{L^1}\quad \text{ for all }\varphi\in\mathcal{D}_0.
\]
Then $T1\in L^\infty$ follows from the Hahn-Banach theorem and $(L^1)^*=L^\infty$.
 \end{proof}

 \subsection{Pointwise Sobolev estimates}
 \begin{proof}[Proof of Corollary \ref{lerner thm}]
   By Theorem \ref{1.1}, if $T1\in L^\infty$, then 
   \[
   \abs{\ave{Tf_1,f_2}}\lesssim\sum_{Q\in\mathcal{S}}\abs{Q}\Omega(f_1;Q)\ave{\abs{f_2}}_Q.
   \]
   From the $L^1$ Poincar\'e inequality, i.e., $\Omega(f_1;Q)\lesssim\abs{Q}^{1/d-1}\int_Q\abs{\nabla f_1}$, we have 
   \[
   \abs{\ave{Tf_1,f_2}}\lesssim\iint \abs{f_2(x)}\abs{\nabla f_1(y)}\sum_{\substack{Q\in\mathcal{S}\\ x,y\in Q}}\ell(Q)^{1-d}\ud y\ud x.
   \]
   Note that if $x,y\in Q$, then it must be $\ell(Q)\gtrsim \abs{x-y}$. Since $\mathcal{S}\subset\mathcal{D}^\sigma$ for some fixed $\sigma$, at each scale there is at most one cube containing a fixed point. Hence for any $x,y$, we can write 
   \begin{align*}
      \sum_{\substack{Q\in\mathcal{S}\\ x,y\in Q}}\ell(Q)^{1-d}=\sum_{k:2^k\gtrsim \abs{x-y} }\sum_{\substack{Q\in\mathcal{S},\ell(Q)=2^k\\ x,y\in Q}}2^{(1-d)k}\lesssim\abs{x-y}^{1-d}.
   \end{align*}
So we obtain $\abs{\ave{Tf_1,f_2}}\lesssim \int I_1(\abs{\nabla f_1})\abs{f_2}.$

 Fix a nonnegative $\rho\in C_c^\infty(\R^d)$ with $\int\rho=1$, and set $\rho_\varepsilon(y)=\varepsilon^{-d}\rho(y/\varepsilon)$. Testing with $f_2(y)=\rho_\varepsilon(x-y)$ gives
\[
\abs{(Tf_1)*\rho_\varepsilon(x)}
\lesssim
\big(I_1(\abs{\nabla f_1})*\rho_\varepsilon\big)(x).
\]
Since both functions are locally integrable, letting \(\varepsilon\to0\) at their common Lebesgue points proves the desired inequality.
 \end{proof}
 \subsection*{Acknowledgements}
This paper is part of the author's PhD thesis. The author is grateful to his doctoral advisor Kangwei Li for introducing him to this problem and for numerous valuable discussions throughout the preparation of this paper. 
 \section*{Declaration of AI use}
The author used ChatGPT to assist with English editing and with the
mathematical development of Section $4$. In particular, ChatGPT help to
explore applications of the oscillatory sparse bounds and refine parts of the proofs. The author takes full
responsibility for the mathematical content and the final manuscript.
\bibliographystyle{plain}
\bibliography{reff}
\end{document}